\documentclass[11pt,letterpaper]{article}
\RequirePackage[backend=biber,bibstyle=authoryear,citestyle=authoryear,
natbib=true,sorting=nyt,sortcites=true,maxnames=10,minnames=10,
hyperref=true,abbreviate=false,date=long,isbn=true,eprint=true,
uniquename=init,giveninits=true]{biblatex}
\usepackage{hyphenat,graphicx}
\usepackage[colorlinks,linktocpage,breaklinks]{hyperref}
\makeatletter
\def\@filecolor{blue}
\def\@linkcolor{blue}
\def\@citecolor{blue}
\def\@urlcolor{blue}
\let\@old@citep\citep
\let\@old@citet\citet
\let\@old@citeauthor\citeauthor
\def\citep{\@old@citep*}
\def\citet{\@old@citet*}
\def\citeauthor{\@old@citeauthor*}
\let\cite\citep
\makeatother
\usepackage[top=1.5in,bottom=1.2in,left=1.35in,right=1.35in,letterpaper]%
  {geometry}
\makeatletter
\def\ps@headings{%
    \let\@mkboth\@gobbletwo
    \def\@oddhead{\hss\scshape\shorttitle\hss\reset@font\rmfamily\thepage}
    \def\@evenhead{\reset@font\rmfamily\thepage\hss\scshape\shortauthors\hss}
    \let\@oddfoot\@empty\let\@evenfoot\@empty}
\@twosidetrue
\makeatother
\usepackage[neveradjust]{paralist}
\usepackage{theorem}
\newtheorem{theorem}{Theorem}[section]
\newtheorem{proposition}[theorem]{Proposition}
\newtheorem{lemma}[theorem]{Lemma}
\newtheorem{prooflemma}{Lemma}[theorem]
{\theorembodyfont{\normalfont\rmfamily}
\newtheorem{definition}[theorem]{Definition}
\newtheorem{example}[theorem]{Example}}
\makeatletter
\newcommand\pushright{\protect\@ADL@pushright}
\newcommand\@ADL@pushright[1]{{\ifvmode\null\hfill{#1}\par\else\ifmmode%
  \@ADLmaths@pushright{\hbox{#1}}\else\ifinner\@ADLhbox@pushright{#1}%
  \else\@ADLparag@pushright{#1}\fi\fi\fi}}
\newcommand\@ADLmaths@pushright[1]{{\ifinner\@ADLhbox@pushright{#1}\else%
  \tag*{$#1$}\fi}}
\newcommand\@ADLparag@pushright[1]{{\parfillskip=0pt\widowpenalty=10000%
  \displaywidowpenalty=10000\finalhyphendemerits=0\@ADLhbox@pushright#1\par}}
\newcommand\@ADLhbox@pushright{\unskip\nobreak\hfil\penalty50\hskip.2em%
  \null\hfill\hfill}
\newenvironment{proof}{\trivlist\item[\hskip\labelsep\textit{Proof:}\/]%
  \@ADLsave@set@qed\xspace\normalfont\rmfamily}
  {\qed\@ADLrestore@qed\endtrivlist}
\newenvironment{subproof}{\trivlist\item[\hskip\labelsep\textit{Proof:}\/]%
  \@ADLsave@set@subqed\normalfont\rmfamily}
  {\subqed\@ADLrestore@subqed\endtrivlist}
\newif\if@ADL@qed\@ADL@qedfalse
\newcommand\qed{\protect\@ADL@qed{$\blacksquare$}}
\newcommand\@ADL@qed[1]{\if@ADL@qed\global\@ADL@qedfalse%
  \pushright{#1}\else\ifhmode\ifinner\else\par\fi\fi\fi}
\newcommand\@ADLrestore@qed{\global\let\if@ADL@qed\@ADLsaved@ifqed}
\newcommand\@ADLsave@set@qed{\let\@ADLsaved@ifqed
  \if@ADL@qed\global\@ADL@qedtrue}

\newif\if@ADL@subqed\@ADL@subqedfalse
\newcommand\subqed{\protect\@ADL@subqed{$\blacktriangledown$}}
\newcommand\@ADL@subqed[1]{\if@ADL@subqed\global\@ADL@subqedfalse%
  \pushright{#1}\else\ifhmode\ifinner\else\par\fi\fi\fi}
\newcommand\@ADLrestore@subqed{\global\let\if@ADL@subqed\@ADLsaved@ifsubqed}
\newcommand\@ADLsave@set@subqed{\let\@ADLsaved@ifsubqed
  \if@ADL@subqed\global\@ADL@subqedtrue}
  
\makeatother
\makeatletter
\newif\if@ADL@oprocend\@ADL@oprocendfalse
\newcommand\oprocend{\@ADLsave@set@oprocend
  \protect\@ADL@oprocend{$\bullet$}\@ADLrestore@oprocend}
\newcommand\@ADL@oprocend[1]{\if@ADL@oprocend\global\@ADL@oprocendfalse%
  \pushright{#1}\else\ifhmode\ifinner\else\par\fi\fi\fi}
\newcommand\@ADLrestore@oprocend{\global
  \let\if@ADL@oprocend\@ADLsaved@ifoprocend}
\newcommand\@ADLsave@set@oprocend{\let\@ADLsaved@ifoprocend\if@ADL@oprocend%
  \global\@ADL@oprocendtrue}
\newcommand\eqoprocend{\tag*{\oprocend}}
\makeatother
\usepackage[reqno,centertags]{amsmath}
\usepackage{amssymb,xspace,mathrsfs}
\let\subset\subseteq

\let\supset\supseteq
\makeatletter
\def\interval{\@ifnextchar({\@ADL@openleftint}{\@ADL@closedleftint}}
\def\@ADL@openleftint(#1,#2{(#1,#2%
  \@ifnextchar){\@ADL@openrightint}{\@ADL@closedrightint}}
\def\@ADL@closedleftint[#1,#2{[#1,#2%
  \@ifnextchar){\@ADL@openrightint}{\@ADL@closedrightint}}
\def\@ADL@openrightint){)}
\def\@ADL@closedrightint]{]}
\makeatother
\newenvironment{keywords}{\quote\small\textbf{Keywords.}}{\endquote}
\newenvironment{AMS}{\quote\small\textbf{AMS Subject Classifications (2010).}}
   {\endquote}
\newcommand\defn[1]{{\normalfont\bfseries\emph{\mathversion{bold}#1}}}
\newcommand\card{\mathupper{card}}
\newcommand\integer{\mathbb{Z}}
\newcommand\integerp{\integer_{>0}}
\newcommand\rational{\mathbb{Q}}
\newcommand\real{\mathbb{R}}
\newcommand\realp{\real_{>0}}
\newcommand\realnn{\real_{\ge0}}
\newcommand\complex{\mathbb{C}}
\newcommand\C{\mathsf{C}}
\renewcommand\L{\mathsf{L}}
\newcommand\M{\mathsf{M}}
\newcommand\R{\mathsf{R}}
\newcommand\lebmes{\lambda}
\newcommand\alg[1]{\mathsf{#1}}
\newcommand\ie{i.e.,}
\newcommand\eg{e.g.,}
\newcommand\cf{cf.}

\newcommand\resp{resp.}
\newcommand\mathupper[1]{\textup{#1}}
\newcommand\eul{\mathupper{e}}
\newcommand\imag{\mathupper{i}}
\newcommand\subscr[2]{#1_{\textup{#2}}}
\renewcommand\d[1]{{\normalfont\textrm{d}}#1}
\newcommand\boundary{\operatorname{boundary}}
\newcommand\closure{\operatorname{closure}}
\newcommand\image{\operatorname{image}}
\newcommand\scirc{\raise1pt\hbox{$\,\scriptstyle\circ\,$}}
\newcommand\map[3]{#1\colon#2\rightarrow#3}

\newcommand\slnorm{\lvert}
\newcommand\srnorm{\rvert}
\newcommand\snorm[1]{\slnorm #1\srnorm}
\newcommand\asnorm[1]{\left\slnorm #1\right\srnorm}
\newcommand\dlnorm{\lVert}
\newcommand\drnorm{\rVert}
\newcommand\dnorm[1]{\dlnorm #1\drnorm}

\newcommand\setdef[2]{\{#1\;|\enspace#2\}}
\newcommand\asetdef[2]{\left\{#1\immediate\vphantom{#2}\;\right|
  \left.\immediate\vphantom{#1}\enspace#2\right\}}

\newcommand\ifam[1]{(#1)}  
\newcommand\pldblref[2]{\mbox{\ref{#1}(\ref{#2})}}  
\title{Should we fly in the Lebesgue-designed airplane?\textemdash{}The
correct defence of the Lebesgue integral\thanks{Research supported in part by
a grant from the Natural Sciences and Engineering Research Council of
Canada}}
\author{Andrew D.\ Lewis\thanks{Professor, Department of Mathematics and
Statistics, Queen's University, Kingston, ON K7L 3N6, Canada,
email:~\texttt{andrew.lewis@queensu.ca}}}
\newcommand\shorttitle{Should we fly in the Lebesgue-designed airplane?}
\newcommand\shortauthors{A.\ D.\ Lewis}
\date{2008/03/04}
\begin{document}
\maketitle

\begin{abstract}
It is well-known that the Lebesgue integral generalises the Riemann integral.
However, as is also well-known but less frequently well-explained, this
generalisation alone is not the reason why the Lebesgue integral is important
and needs to be a part of the arsenal of any mathematician, pure or applied.
Those who understand the correct reasons for the importance of the Lebesgue
integral realise there are at least two crucial differences between the
Riemann and Lebesgue theories.  One is the difference between the Dominated
Convergence Theorem in the two theories, and another is the completeness of
the normed vector spaces of integrable functions.  Here topological
interpretations are provided for the differences in the Dominated Convergence
Theorems, and explicit counterexamples are given which illustrate the
deficiencies of the Riemann integral.  Also illustrated are the deleterious
consequences of the defects in the Riemann integral on Fourier transform
theory if one restricts to Riemann integrable functions.
\end{abstract}
\begin{keywords}
Lebesgue integral, Riemann integral.
\end{keywords}
\begin{AMS}
28-01
\end{AMS}

\section{Introduction}

The title of this paper is a reference to the well-known quote of the applied
mathematician and engineer Richard W.~Hamming (1915--1998):
\begin{quote}
Does anyone believe that the difference between the Lebesgue and Riemann
integrals can have physical significance, and that whether say, an airplane
would or would not fly could depend on this difference? If such were claimed,
I should not care to fly in that plane.
\end{quote}
The statement by \citeauthor{RWH:80} is open to many interpretations, but the
interpretation of \citeauthor{RWH:80} himself can be gleaned
from~\cite{RWH:80} and, particularly,~\cite{RWH:98}\@; also see~\cite{PJD:98}
for a discussion of some of \citeauthor{RWH:80}'s views on mathematics and
the ``real world.''  Perhaps a fair summary of \citeauthor{RWH:98}'s views
toward the Riemann and Lebesgue theories of integration is that the
distinction between them is not apt to be seen in Nature.  This seems about
right to us.  Unfortunately, however, this quote of \citeauthor{RWH:98}'s is
often used in a confused manner that indicates the quoter's misunderstanding
of the purpose and importance of the Lebesgue integral.  Indeed, very often
\citeauthor{RWH:98}'s quote is brought out as an excuse to disregard Lebesgue
integration, the idea being that it is the product of some vapid pursuit of
generality.  Oxymoronically, this is often done simultaneously with the free
use of results (like completeness of the $\L^p$-spaces) which rely crucially
on the distinctions between the Riemann and Lebesgue integrals.

That the value of the Lebesgue theory of integration (and all of the theories
equivalent to or generalising it\footnote{\citet{RWH:98} himself seems to
find the Henstock integral, which generalises the Lebesgue integral, to be a
more satisfactory construction than the Lebesgue integral.}) may not be
appreciated fully by non-mathematicians should not be too surprising: the
Lebesgue theory is subtle.  Moreover, it is definitely not the case that the
importance of the Lebesgue theory over the Riemann theory is explained
clearly in all texts on integration theory; in fact, the important
distinctions are rarely explicitly stated, though they are almost always
implicitly present.  What is most discomforting, however, is that
mathematicians themselves sometimes offer an \emph{incorrect} defence of the
Lebesgue theory.  For example, it is not uncommon to see defences made along
the lines of, ``The class of functions that can be integrated using the
Lebesgue theory is larger than that using the Riemann
theory''~\cite{MD/MI:02}\@.  Sometimes, playing into the existing suspicions
towards unnecessary generality, it is asserted that the mere fact that the
Lebesgue theory generalises the Riemann theory is sufficient to explain its
superiority.  These sorts of defences of the Lebesgue theory are certainly
factual.  But they are also emphatically \emph{not} the sorts of reasons why
the Lebesgue theory is important.  The functions that can be integrated using
the Lebesgue theory, but which cannot be integrated using the Riemann theory,
are not important; just try showing a Lebesgue integrable but not Riemann
integrable function to someone who is interested in applications of
mathematics, and see if they think the function is important.  This certainly
must at least partially underlie \citeauthor{RWH:98}'s motivation for his
quote.

The value of the Lebesgue theory over the Riemann theory is that it is
superior, as a \emph{theory of integration}\@.  By this it is meant that
there are theorems in the Lebesgue theory that are true and useful, but that
are not true in the Riemann theory.  Probably the most crucial such theorem
is the powerful version of the Dominated Convergence Theorem that one has in
the Lebesgue theory.  This theorem is constantly used in the proof of many
results that are important in applications.  For example, the Dominated
Convergence Theorem is used crucially in the proof of the completeness of
$\L^p$-spaces.  In turn, the completeness of these spaces is an essential
part of why these spaces are useful in, for example, the theory of signals
and systems that is taught to all Electrical Engineering undergraduates.  For
instance, in many texts on signals and systems one can find the statement
that the Fourier transform is an isomorphism of $\L^2(\real;\complex)$\@.
This statement is one that needs a sort of justification that is
(understandably) not often present in such texts.  But its presence can at
least be justified by its being correct.  With only the Riemann theory of
integration at one's disposal, the statement is simply not correct.  We
illustrate this in Section~\ref{subsec:riemann-ccft}\@.

In this paper we provide topological explanations for the differences between the Riemann and Lebesgue theories of integration.  The intent is not to make these differences clear for non-mathematicians.  Indeed, for non-mathematicians the contents of the preceding paragraphs, along with the statement that, ``The Lebesgue theory of integration is to the Riemann theory of integration as the real numbers are to the rational numbers,'' (this is the content of our Example~\ref{eg:Riemann!Cauchy}) seems about the best one can do.  No, what we aim to do in this paper is clarify \emph{for mathematicians} the reasons for the superiority of the Lebesgue theory.  We do this by providing two theorems, both topological in nature, that are valid for the Lebesgue theory and providing counterexamples illustrating that they are not true for the Riemann theory.  We also illustrate the consequences of the topological deficiencies of the Riemann integral by explicitly depicting the limitations of the $\L^2$-Fourier transform with the Riemann integral.  One of the contributions in this paper is that we give the \emph{correct} counterexamples.  All too often one sees counterexamples that illustrate \emph{some} point, but not always the one that one wishes to make.  The core of what we say here exists in the literature in one form or another, although we have not seen in the literature counterexamples that illustrate what we illustrate in Examples~\ref{eg:R1!hatMp-closed} and~\ref{eg:Riemann!Cauchy}\@.  The principal objective here is to organise the results and examples in an explicit and compelling way.

In the event that the reader is consulting this paper in a panic just prior
to boarding an airplane, let us answer the question posed in the title of the
paper.  The answer is, ``The question is meaningless as the distinctions
between the Riemann and Lebesgue integrals do not, and should not be thought
to, contribute to such worldly matters as aircraft design.''  However, the
salient point is that this is not a valid criticism of the Lebesgue integral.
What follows is, we hope, a valid defence of the Lebesgue integral.

\section{Spaces of functions}\label{sec:LR-topologies}

To keep things simple and to highlight the important parts of our
presentation, in this section and in most of the rest of the paper we
consider $\real$-valued functions defined on
$I=\interval[0,1]\subset\real$\@.  Extensions to more general settings are
performed mostly by simple changes of notation.  The Lebesgue measure on
$\interval[0,1]$ is denoted by $\lebmes$\@.  In order to distinguish the
Riemann and Lebesgue integrals we denote them by
\begin{equation*}
\int_0^1f(x)\,\d{x},\qquad\int_If\,\d{\lebmes},
\end{equation*}
respectively.  In order to make our statements as strong as possible, by the
Riemann integral we mean the improper Riemann integral to allow for unbounded
functions~\cite[Section~8.5]{JEM/MJH:93}\@.

\subsection{Normed vector spaces of integrable
functions}\label{subsec:nvs-integrable}

We use slightly unconventional notation that is internally self-consistent
and convenient for our purposes here.

Let us first provide the large vector spaces whose subspaces are of interest.
By $\real^{\interval[0,1]}$ we denote the set of all $\real$-valued functions
on $\interval[0,1]$\@.  This is also the product of $\card(\interval[0,1])$
copies of $\real$\@, and we shall alternatively think of elements of
$\real^{\interval[0,1]}$ as being functions or elements of a product of sets,
as is convenient for our purposes.  We consider the standard $\real$-vector
space structure on $\real^{\interval[0,1]}$\@:
\begin{equation*}
(f+g)(x)=f(x)+g(x),\quad(af)(x)=a(f(x)),\qquad f,g\in\real^{\interval[0,1]},\
a\in\real.
\end{equation*}

In $\real^{\interval[0,1]}$ consider the subspace
\begin{equation*}
Z(\interval[0,1];\real)=\asetdef{f\in\real^{\interval[0,1]}}
{\lambda(\setdef{x\in\interval[0,1]}{f(x)\not=0})=0}.
\end{equation*}
Then denote $\hat{\real}^{\interval[0,1]}=\real^{\interval[0,1]}/
Z(\interval[0,1];\real)$\@; this is then the set of equivalence classes of
functions agreeing almost everywhere.  We shall denote by
$[f]=f+Z(\interval[0,1];\real)$ the equivalence class of
$f\in\real^{\interval[0,1]}$\@.

Now let us consider subspaces of $\real^{\interval[0,1]}$ and
$\hat{\real}^{\interval[0,1]}$ consisting of integrable functions.  Let us
denote by
\begin{equation*}
\R^1(\interval[0,1];\real)=\setdef{\map{f}{\interval[0,1]}{\real}}
{f\ \textrm{is Riemann integrable}}.
\end{equation*}
On $\R^1(\interval[0,1];\real)$ define a seminorm $\dnorm{\cdot}_1$ by
\begin{equation*}
\dnorm{f}_1=\int_0^1\snorm{f(x)}\,\d{x},
\end{equation*}
and denote
\begin{equation*}
\R_0(\interval[0,1];\real)=\asetdef{f\in\R^1(\interval[0,1];\real)}
{\dnorm{f}_1=0}.
\end{equation*}
Then we define
\begin{equation*}
\hat{\R}^1(\interval[0,1];\real)=\R^1(\interval[0,1];\real)/
\R_0(\interval[0,1];\real),
\end{equation*}
and note that this $\real$-vector space is then equipped with the norm
$\dnorm{[f]}_1=\dnorm{f}_1$\@, accepting the slight abuse of notation of
using $\dnorm{\cdot}_1$ is different contexts.

The preceding constructions can be carried out, replacing ``$\R$'' with
``$\L$'' and replacing the Riemann integral with the Lebesgue integral, to
arrive at the seminormed vector space
\begin{equation*}
\L^1(\interval[0,1];\real)=\setdef{\map{f}{\interval[0,1]}{\real}}
{f\ \textrm{is Lebesgue integrable}}
\end{equation*}
with the seminorm
\begin{equation*}
\dnorm{f}_1=\int_I\snorm{f}\,\d{\lebmes},
\end{equation*}
the subspace
\begin{equation*}
\L_0(\interval[0,1];\real)=\asetdef{f\in\L^1(\interval[0,1];\real)}
{\dnorm{f}_1=0},
\end{equation*}
and the normed vector space
\begin{equation*}
\hat{\L}^1(\interval[0,1];\real)=\L^1(\interval[0,1];\real)/
\L_0(\interval[0,1];\real).
\end{equation*}
We denote the norm on $\hat{\L}^1(\interval[0,1];\real)$ by
$\dnorm{\cdot}_1$\@, this not being too serious an abuse of notation since
$\hat{\R}^1(\interval[0,1];\real)$ is a subspace of
$\hat{\L}^1(\interval[0,1];\real)$ with the restriction of the norm on
$\hat{\L}^1(\interval[0,1];\real)$ to $\hat{\R}^1(\interval[0,1];\real)$
agreeing with the norm on $\hat{\R}^1(\interval[0,1];\real)$\@.  This is a
consequence of the well-known fact that the Lebesgue integral generalises the
Riemann integral~\cite[Theorem~2.5.4]{DLC:13}\@.

During the course of the development of the Lebesgue theory of integration
one shows that
\begin{equation*}
\L_0(\interval[0,1];\real)=Z(\interval[0,1];\real)
\end{equation*}
\cite[\eg][Corollary~2.3.11]{DLC:13}\@.  The corresponding assertion is not
true for the Riemann theory.
\begin{example}\label{eg:Qcharfunc}
Let us denote by $F\in\real^{\interval[0,1]}$ the characteristic function of
$\rational\cap\interval[0,1]$\@.  This is perhaps the simplest and most
commonly seen example of a function that is Lebesgue integrable but not
Riemann integrable~\cite[Example~2.5.4]{DLC:13}\@.  Thus
$F\not\in\R_0(\interval[0,1];\real)$\@.  However, $F\in
Z(\interval[0,1];\real)$\@.\oprocend
\end{example}

While the preceding example is often used as an example of a function that is
not Riemann integrable but is Lebesgue integrable, one needs to be careful
about exaggerating the importance, even mathematically, of this example.  In
Examples~\ref{eg:R1!hatMp-closed} and~\ref{eg:Riemann!Cauchy} below we
shall see that this example is not sufficient for demonstrating some of the
more important deficiencies of the Riemann integral.

\subsection{Pointwise convergence topologies}

For $x\in\interval[0,1]$ let us denote by
$\map{p_x}{\real^{\interval[0,1]}}{\realnn}$ the seminorm defined by
$p_x(f)=\snorm{f(x)}$\@.  The family $\ifam{p_x}_{x\in\interval[0,1]}$ of
seminorms on $\real^{\interval[0,1]}$ defines a locally convex topology.  A
basis of open sets for this topology is given by products of the form
$\prod_{x\in\interval[0,1]}U_x$ where $U_x\subset\real$ is open and where
$U_x=\real$ for all but finitely many $x\in\interval[0,1]$\@.  A sequence
$\ifam{f_j}_{j\in\integerp}$ in $\real^{\interval[0,1]}$ converges to
$f\in\real^{\interval[0,1]}$ if and only if the sequence converges pointwise
in the usual sense~\cite[Theorem~42.2]{SW:70}\@.  Let us, therefore, call
this the \defn{topology of pointwise convergence} and let us denote by
$\C_p(\interval[0,1];\real)$ the vector space $\real^{\interval[0,1]}$ when
equipped with this topology.  For clarity, we shall prefix with ``$\C_p$''
topological properties in the topology of pointwise convergence.  Thus, for
example, we shall say ``$\C_p$-open'' to denote an open set in
$\C_p(\interval[0,1];\real)$\@.

We will be interested in bounded subsets of $\C_p(\interval[0,1];\real)$\@.
We shall use the characterisation of a bounded subset $B$ of a topological
$\real$-vector space $\alg{V}$ that a set is bounded if and only if, for
every sequence $\ifam{v_j}_{j\in\integerp}$ in $B$ and for every sequence
$\ifam{a_j}_{j\in\integerp}$ in $\real$ converging to $0$\@, it holds that
the sequence $\ifam{a_jv_j}_{j\in\integerp}$ converges to zero in the
topology of $\alg{V}$~\cite[Theorem~1.30]{WR:91}\@.
\begin{proposition}\label{prop:Cp-bounded}
A subset\/ $B\subset\C_p(\interval[0,1];\real)$ is\/ $\C_p$-bounded if and
only if there exists a nonnegative-valued\/ $g\in\real^{\interval[0,1]}$ such
that
\begin{equation*}
B\subset\setdef{f\in\C_p(\interval[0,1];\real)}
{\snorm{f(x)}\le g(x)\ \textrm{for every}\ x\in\interval[0,1]}.
\end{equation*}
\begin{proof}
Suppose that there exists a nonnegative-valued $g\in\real^{\interval[0,1]}$
such that $\snorm{f(x)}\le g(x)$ for every $x\in\interval[0,1]$ if $f\in
B$\@.  Let $\ifam{f_j}_{j\in\integerp}$ be a sequence in $B$ and let
$\ifam{a_j}_{j\in\integerp}$ be a sequence in $\real$ converging to $0$\@.
If $x\in\interval[0,1]$ then
\begin{equation*}
\lim_{j\to\infty}\snorm{a_jf_j(x)}\le\lim_{j\to\infty}\snorm{a_j}g(x)=0,
\end{equation*}
which gives $\C_p$-convergence of the sequence
$\ifam{a_jf_j}_{j\in\integerp}$ to zero.

Next suppose that there exists no nonnegative-valued function
$g\in\real^{\interval[0,1]}$ such that $\snorm{f(x)}\le g(x)$ for every
$x\in\interval[0,1]$ if $f\in B$\@.  This means that there exists
$x_0\in\interval[0,1]$ such that, for every $M\in\realp$\@, there exists
$f\in B$ such that $\snorm{f(x_0)}>M$\@.  Let $\ifam{a_j}_{j\in\integerp}$ be
a sequence in $\real$ converging to $0$ and such that $a_j\not=0$ for every
$j\in\integerp$\@.  Then let $\ifam{f_j}_{j\in\integerp}$ be a sequence in
$B$ such that $\snorm{f_j(x_0)}>\asnorm{a_j^{-1}}$ for every
$j\in\integerp$\@.  Then $\snorm{a_jf_j(x_0)}>1$ for every $j\in\integerp$\@,
implying that the sequence $\ifam{a_jf_j}_{j\in\integerp}$ cannot
$\C_p$-converge to zero.  Thus $B$ is not $\C_p$-bounded.
\end{proof}
\end{proposition}

Of course, in the theory of integration one is interested, not in pointwise
convergence of arbitrary functions, but in pointwise convergence of
measurable functions.  Let us, therefore, denote
\begin{equation*}
\M(\interval[0,1];\real)=\asetdef{f\in\real^{\interval[0,1]}}
{f\ \textrm{is Lebesgue measurable}},
\end{equation*}
where we understand the topology on $\M(\interval[0,1];\real)$ to be the
subspace topology inherited from $\C_p(\interval[0,1];\real)$\@.  Standard
theorems on measurable functions show that $\M(\interval[0,1];\real)$ is a
subspace~\cite[Proposition~2.1.6]{DLC:13} and is $\C_p$-sequentially
closed~\cite[Proposition~2.1.5]{DLC:13}.  Moreover, the stronger assertion of
closedness does not hold.  The following result shows this, as well as giving
topological properties of $Z(\interval[0,1];\real)$\@.
\begin{proposition}
The subspaces\/ $\M(\interval[0,1];\real)$ and\/ $Z(\interval[0,1];\real)$
of\/ $\C_p(\interval[0,1];\real)$ are not\/ $\C_p$-closed, but are\/
$\C_p$-sequentially closed.
\begin{proof}
The $\C_p$-sequential closedness of $\M(\interval[0,1];\real)$ and
$Z(\interval[0,1];\real)$ follow from standard theorems, as pointed out
above.  We first show that
$\C_p(\interval[0,1];\real)\setminus\M(\interval[0,1];\real)$ is not
$\C_p$-open.  Let
$f\in\C_p(\interval[0,1];\real)\setminus\M(\interval[0,1];\real)$ and let $V$
be a $\C_p$-open set containing $f$\@.  Then $V$ contains a basic
neighbourhood.  Thus there exists $\epsilon\in\realp$\@,
$x_1,\dots,x_k\in\interval[0,1]$\@, and a basic neighbourhood
$U=\prod_{x\in\interval[0,1]}U_x$ contained in $V$ where
\begin{enumerate}
\item $U_{x_j}=\interval({f(x_j)-\epsilon},{f(x_j)+\epsilon})$ for
$j\in\{1,\dots,k\}$ and
\item $U_x=\real$ for $x\in\interval[0,1]\setminus\{x_1,\dots,x_k\}$\@.
\end{enumerate}
Then the function $g\in\real^{\interval[0,1]}$ defined by
\begin{equation*}
g(x)=\begin{cases}f(x),&x\in\{x_1,\dots,x_k\},\\
0,&\textrm{otherwise}\end{cases}
\end{equation*}
is in $U\cap\M(\interval[0,1];\real)$\@.  Thus $g\in V$\@, so showing that
every neighbourhood of $f$ contains a member of
$\M(\interval[0,1];\real)$\@.

To show that $Z(\interval[0,1];\real)$ is not $\C_p$-closed we shall show
that $\C_p(\interval[0,1];\real)\setminus Z(\interval[0,1];\real)$ is not
$\C_p$-open.  Let $f\in\C_p(\interval[0,1];\real)\setminus
Z(\interval[0,1];\real)$ be given by $f(x)=1$ for all $x\in\interval[0,1]$\@.
Let $V$ be a $\C_p$-open subset containing $f$\@.  Then $V$ contains a basic
neighbourhood from $\C_p(\interval[0,1];\real)$\@, and in particular a basic
neighbourhood of the form $U=\prod_{x\in\interval[0,1]}U_x$ where the
$\C_p$-open sets $U_x\subset\real$\@, $x\in\interval[0,1]$\@, have the
following properties:
\begin{enumerate}
\item there exists $\epsilon\in\interval(0,1)$ and a finite set
$x_1,\dots,x_k\in\interval[0,1]$ such that
$U_{x_j}=\interval({1-\epsilon},{1+\epsilon})$ for each
$j\in\{1,\dots,k\}$\@;
\item for $x\in\interval[0,1]\setminus\{x_1,\dots,x_k\}$ we have
$U_x=\real$\@.
\end{enumerate}
We claim that such a basic neighbourhood $U$ contains a function from
$Z(\interval[0,1];\real)$\@.  Indeed, the function
\begin{equation*}
g(x)=\begin{cases}1,&x\in\{x_1,\dots,x_k\},\\0,&\textrm{otherwise}\end{cases}
\end{equation*}
is in $U\cap Z(\interval[0,1];\real)$\@, and so is in $V\cap
Z(\interval[0,1];\real)$\@.  This shows that
$\C_p(\interval[0,1];\real)\setminus Z(\interval[0,1];\real)$ is not
$\C_p$-open, as desired.
\end{proof}
\end{proposition}

\subsection{Almost everywhere pointwise convergence limit structures}

For many applications, it is the space $\hat{\L}^1(\interval[0,1];\real)$\@,
not $\L(\interval[0,1];\real)$\@, that is of interest, this by virtue of its
possessing a norm and not a seminorm.  (Of course, one might also be
interested in $\hat{\R}^1(\interval[0,1];\real)$\@, but the point of this
paper is to clarify the ways in which this space is deficient.)  Thus one is
interested in subspaces of $\hat{\real}^{\interval[0,1]}$\@.  The largest
such subspace in which we shall be interested is the image of the Lebesgue
measurable functions in $\hat{\real}^{\interval[0,1]}$ under the projection
from $\real^{\interval[0,1]}$\@:
\begin{equation*}
\hat{\M}(\interval[0,1];\real)=
\M(\interval[0,1];\real)/Z(\interval[0,1];\real).
\end{equation*}
Note that the quotient is well-defined since completeness of the Lebesgue
measure gives $Z(\interval[0,1];\real)\subset\M(\interval[0,1];\real)$\@.

Now, one wishes to provide structure on $\hat{\M}(\interval[0,1];\real)$
such that there is a notion of convergence which agrees with almost
everywhere pointwise convergence.  First let us be clear about what we mean
by almost everywhere pointwise convergence relative to the various function
spaces we are using.
\begin{definition}
\begin{compactenum}[(i)]
\item A sequence $\ifam{f_j}_{j\in\integerp}$ in $\M(\interval[0,1];\real)$
is \defn{almost everywhere pointwise convergent} to
$f\in\M(\interval[0,1];\real)$ if
\begin{equation*}
\lambda(\setdef{x\in\interval[0,1]}{\ifam{f_j(x)}\ \textrm{does not converge
to}\ f(x)})=0.
\end{equation*}
\item A sequence $\ifam{[f_j]}_{j\in\integerp}$ in
$\hat{\M}(\interval[0,1];\real)$ is \defn{almost everywhere pointwise
convergent} to $[f]\in\hat{\M}(\interval[0,1];\real)$ if
\begin{equation*}\eqoprocend
\lambda(\setdef{x\in\interval[0,1]}{\ifam{f_j(x)}\ \textrm{does not converge
to}\ f(x)})=0.
\end{equation*}
\end{compactenum}
\end{definition}

We should ensure that the definition of almost everywhere pointwise
convergence in $\hat{\M}(\interval[0,1];\real)$ is well-defined.
\begin{lemma}\label{lem:pwaec}
For a sequence\/ $\ifam{[f_j]}_{j\in\integerp}$ in\/
$\hat{\M}(\interval[0,1];\real)$ and for\/
$[f]\in\hat{\M}(\interval[0,1];\real)$ the following statements are
equivalent:
\begin{compactenum}[(i)]
\item there exists a sequence\/ $\ifam{g_j}_{j\in\integerp}$ in\/
$\M(\interval[0,1];\real)$ and\/ $g\in\M(\interval[0,1];\real)$ such that
\begin{compactenum}[(a)]
\item $[g_j]=[f_j]$ for\/ $j\in\integerp$\@,
\item $[g]=[f]$\@, and
\item $\ifam{g_j}_{j\in\integerp}$ converges pointwise almost everywhere to\/
$g$\@.
\end{compactenum}
\item for every sequence\/ $\ifam{g_j}_{j\in\integerp}$ in\/
$\M(\interval[0,1];\real)$ and for every\/ $g\in\M(\interval[0,1];\real)$
satisfying
\begin{compactenum}[(a)]
\item $[g_j]=[f_j]$ for\/ $j\in\integerp$ and
\item $[g]=[f]$\@,
\end{compactenum}
it holds that\/ $\ifam{g_j}_{j\in\integerp}$ converges pointwise almost
everywhere to\/ $g$\@.
\end{compactenum}
\begin{proof}
It is clear that the second statement implies the first, so we only prove the
converse.  Thus we let $\ifam{g_j}_{j\in\integerp}$ in
$\M(\interval[0,1];\real)$ and\/ $g\in\M(\interval[0,1];\real)$ be such that
\begin{enumerate}
\item $[g_j]=[f_j]$ for\/ $j\in\integerp$\@,
\item $[g]=[f]$\@, and
\item $\ifam{g_j}_{j\in\integerp}$ converges pointwise almost everywhere to
$g$\@.
\end{enumerate}
Let $\ifam{h_j}_{j\in\integerp}$ be a sequence in $\M(\interval[0,1];\real)$
and let $h\in\M(\interval[0,1];\real)$ be such that
\begin{enumerate}
\item $[h_j]=[f_j]$ for\/ $j\in\integerp$ and
\item $[h]=[f]$\@.
\end{enumerate}
Define
\begin{equation*}
A=\setdef{x\in\interval[0,1]}{g(x)\not=f(x)},\quad
B=\setdef{x\in\interval[0,1]}{h(x)\not=f(x)}
\end{equation*}
and, for $j\in\integerp$\@, define
\begin{equation*}
A_j=\setdef{x\in\interval[0,1]}{g_j(x)\not=f_j(x)},\quad
B_j=\setdef{x\in\interval[0,1]}{h_j(x)\not=f_j(x)}
\end{equation*}
and note that
\begin{equation*}
x\in\interval[0,1]\setminus(A\cup B)=(\interval[0,1]\setminus A)\cap
(\interval[0,1]\setminus B)\quad\implies\quad h(x)=f(x)=g(x)
\end{equation*}
and
\begin{equation*}
x\in\interval[0,1]\setminus(A_j\cup B_j)=(\interval[0,1]\setminus
A_j)\cap(\interval[0,1]\setminus B_j)\quad\implies\quad
h_j(x)=f_j(x)=g_j(x).
\end{equation*}
Thus,
\begin{equation*}
x\in\interval[0,1]\setminus\left((\cup_{j\in\integerp}A_j\cup B_j)\cup
(A\cup B)\right)\quad\implies\quad
\lim_{j\to\infty}h_j(x)=\lim_{j\to\infty}g_j(x)=g(x)=h(x).
\end{equation*}
Since $(\cup_{j\in\integerp}A_j\cup B_j)\cup(A\cup B)$ is a countable union
of sets of measure zero, it has zero measure, and so
$\ifam{h_j}_{j\in\integerp}$ converges pointwise almost everywhere to $h$\@.
\end{proof}
\end{lemma}

Now that we understand just what sort of convergence we seek in
$\hat{\M}(\interval[0,1];\real)$\@, we can think about how to achieve this.
The obvious first guess is to use the quotient topology on
$\hat{\M}(\interval[0,1];\real)$ inherited from the $\C_p$-topology on
$\M(\interval[0,1];\real)$\@.  However, convergence in this topology fails to
agree with almost everywhere pointwise convergence.  Indeed, we have the
following more sweeping statement.
\begin{proposition}\label{prop:aetop!exist}
Let\/ $\subscr{\mathscr{T}}{a.e.}$ be the set of topologies\/ $\tau$ on\/
$\hat{\M}(\interval[0,1];\real)$ such that the convergent sequences in\/
$\tau$ are precisely the almost everywhere pointwise convergent sequences.
Then\/ $\subscr{\mathscr{T}}{a.e.}=\emptyset$\@.
\begin{proof}
Suppose that $\subscr{\mathscr{T}}{a.e.}\not=\emptyset$ and let
$\tau\in\subscr{\mathscr{T}}{a.e.}$\@.  Let us denote by
$z\in\M(\interval[0,1];\real)$ the zero function.  Let
$\ifam{f_j}_{j\in\integerp}$ be a sequence in $\M(\interval[0,1];\real)$
converging in measure to $z$\@,~\ie~for every $\epsilon\in\realp$\@,
\begin{equation*}
\lim_{j\to\infty}\lambda(\setdef{x\in\interval[0,1]}
{\snorm{f_j(x)}\ge\epsilon})=0,
\end{equation*}
but not converging pointwise almost everywhere to
$z$~\cite[Example~3.1.1(b)]{DLC:13}\@.  Since almost everywhere pointwise
convergence agrees with convergence in $\tau$\@, there exists a neighbourhood
$U$ of $[z]$ in $\hat{\M}(\interval[0,1];\real)$ such that the set
\begin{equation*}
\setdef{j\in\integerp}{[f_j]\in U}
\end{equation*}
is finite.  As is well-known~\cite[Proposition~3.1.3]{DLC:13}\@, there exists
a subsequence $\ifam{f_{j_k}}_{k\in\integerp}$ of
$\ifam{f_j}_{j\in\integerp}$ that converges pointwise almost everywhere to
$z$\@.  Thus the sequence $\ifam{[f_{j_k}]}_{k\in\integerp}$ converges
pointwise almost everywhere to $[z]$\@, and so converges to $[z]$ in
$\tau$\@.  Thus, in particular, the set
\begin{equation*}
\setdef{k\in\integerp}{[f_{j_k}]\in U}
\end{equation*}
is infinite, which is a contradiction.
\end{proof}
\end{proposition}

It is moderately well-known that there can be no topology on $\hat{\M}(\interval[0,1];\real)$ which gives rise to almost everywhere pointwise convergence.  For instance, this is observed by \cite{MF:21}\@.  Our proof of Proposition~\ref{prop:aetop!exist} is adapted slightly from the observation of \citet{ETO:66}\@.  The upshot of the result is that, if one is going to provide some structure with which to describe almost everywhere pointwise convergence, this structure must be something different than a topology.  This was addressed by \citet{RA:50} who observed that the notion of convergence in measure \emph{is} topological, but almost everywhere pointwise convergence is not.  To structurally distinguish between the two sorts of convergence, \citeauthor{RA:50} introduces the notion of a limit structure.  This idea is discussed in some generality for Borel measurable functions by \citet{UH:00} using multiple valued topologies.  Here we will introduce the notion of a limit structure in as direct a manner as possible, commensurate with our objectives.  Readers wishing to explore the subject in more detail are referred to \cite{RB/HPB:02}\@.

For a set $X$ let $\mathscr{F}(X)$ denote the set of filters on $X$ and, for
$x\in X$\@, denote by
\begin{equation*}
\mathcal{F}_x=\setdef{S\subset X}{x\in S}
\end{equation*}
the principal filter generated by $\{x\}$\@.  If $(\Lambda,\preceq)$ is a
directed set and if $\map{\phi}{\Lambda}{X}$ is a $\Lambda$-net, we denote
the tails of the net $\phi$ by
\begin{equation*}
T_\phi(\lambda)=\setdef{\phi(\lambda')}
{\lambda'\ge\lambda},\qquad\lambda\in\Lambda.
\end{equation*}
We then denote by
\begin{equation*}
\mathcal{F}_\phi=\setdef{S\subset X}
{T_\phi(\lambda)\subset S\ \textrm{for some}\ \lambda\in\Lambda}
\end{equation*}
the ``tail filter'' (also sometimes called the ``Fr\'echet filter'')
associated to the $\Lambda$-net $\phi$\@.
\begin{definition}
A \defn{limit structure} on a set $X$ is a subset
$\mathscr{L}\subset\mathscr{F}(X)\times X$ with the following properties:
\begin{compactenum}[(i)]
\item \label{pl:ls1} if $x\in X$ then $(\mathcal{F}_x,x)\in\mathscr{L}$\@;
\item \label{pl:ls2} if $(\mathcal{F},x)\in\mathscr{L}$ and if
$\mathcal{F}\subset\mathcal{G}\in\mathscr{F}(X)$ then
$(\mathcal{G},x)\in\mathscr{L}$\@;
\item \label{pl:ls3} if $(\mathcal{F},x),(\mathcal{G},x)\in\mathscr{L}$ then
$(\mathcal{F}\cap\mathcal{G},x)\in\mathscr{L}$\@.
\end{compactenum}
If $(\Lambda,\preceq)$ is a directed set, a $\Lambda$-net
$\map{\phi}{\Lambda}{X}$ is \defn{$\mathscr{L}$-convergent} to $x\in X$ if
$(\mathcal{F}_\phi,x)\in\mathscr{L}$\@.  Let us denote by
$\mathscr{S}(\mathscr{L})$ the set of $\mathscr{L}$-convergent
$\integerp$-nets,~\ie~the set of $\mathscr{L}$-convergent sequences.\oprocend
\end{definition}

The intuition behind the notion of a limit structure is as follows.
Condition~\eqref{pl:ls1} says that the trivial filter converging to $x$
should be included in the limit structure, condition~\eqref{pl:ls2} says that
if a filter converges to $x$\@, then every coarser filter also converges to
$x$\@, and condition~\eqref{pl:ls3} says that ``mixing'' filters converging
to $x$ should give a filter converging to $x$\@. Starting from the definition
of a limit structure, one can reproduce many of the concepts from
topology,~\eg~openness, closedness, compactness, continuity. Since we are not
going to make us of any of these general constructions, we merely refer the
interested reader to \cite{RB/HPB:02}\@. The one notion we will use is the
following: a subset $A$ of a set $X$ with a limit structure $\mathscr{L}$ is
\defn{$\mathscr{L}$-sequentially closed} if every $\mathscr{L}$-convergent
sequence $\ifam{x_j}_{j\in\integerp}$ in $A$ converges to $x\in A$\@.

We are interested in the special case of limit structures on a $\real$-vector
space $\alg{V}$\@; one will trivially see that there is nothing in the
definitions that requires the field to be $\real$\@.  For
$\mathcal{F},\mathcal{G}\in\mathscr{F}(\alg{V})$ and for $a\in\real$ we
denote
\begin{equation*}
\mathcal{F}+\mathcal{G}=\setdef{A+B}{A\in\mathcal{F},\ B\in\mathcal{G}},\quad
a\mathcal{F}=\setdef{aA}{A\in\mathcal{F}},
\end{equation*}
where, as usual,
\begin{equation*}
A+B=\setdef{u+v}{u\in A,\ v\in B},\quad aA=\setdef{au}{u\in A}.
\end{equation*}
We say that a limit structure $\mathscr{L}$ on a vector space $\alg{V}$ is
\defn{linear} if $(\mathcal{F}_1,v_1),(\mathcal{F}_2,v_2)\in\mathscr{L}$
implies that $(\mathcal{F}_1+\mathcal{F}_2,v_1+v_2)\in\mathscr{L}$ and if
$a\in\real$ and $(\mathcal{F},v)\in\mathscr{L}$ then
$(a\mathcal{F},av)\in\mathscr{L}$\@.  Following the characterisation of
bounded subsets of topological vector spaces, we say a subset
$B\subset\alg{V}$ is \defn{$\mathscr{L}$-bounded} if, for every sequence
$\ifam{v_j}_{j\in\integerp}$ in $B$ and for every sequence
$\ifam{a_j}_{j\in\integerp}$ in $\real$ converging to $0$\@, the sequence
$\ifam{a_jv_j}_{j\in\integerp}$ is $\mathscr{L}$-convergent to zero.

For $[f]\in\hat{\M}(\interval[0,1];\real)$ define
\begin{multline*}
\mathscr{F}_{[f]}=\{\mathcal{F}\in
\mathscr{F}(\hat{\M}(\interval[0,1];\real))\mid\enspace
\mathcal{F}_\phi\subset\mathcal{F}\ \textrm{for some}\
\integerp\textrm{-net}\ \phi\ \textrm{such that}\\
\ifam{\phi(j)}_{j\in\integerp}\ \textrm{is almost everywhere pointwise
convergent to}\ [f]\}.
\end{multline*}
We may now define a limit structure on $\hat{\M}(\interval[0,1];\real)$ as
follows.
\begin{theorem}\label{the:hatM-convergence}
The subset of\/ $\mathscr{F}(\hat{\M}(\interval[0,1];\real))\times
\hat{\M}(\interval[0,1];\real)$ defined by
\begin{equation*}
\mathscr{L}_\lambda=
\asetdef{(\mathcal{F},[f])}{\mathcal{F}\in\mathscr{F}_{[f]}}
\end{equation*}
is a linear limit structure on\/ $\hat{\M}(\interval[0,1];\real)$\@.
Moreover, a sequence\/ $\ifam{[f_j]}_{j\in\integerp}$ is\/
$\mathscr{L}_\lambda$-convergent to\/ $[f]$ if and only if the sequence is
almost everywhere pointwise convergent to\/ $[f]$\@.
\begin{proof}
Let $[f]\in\hat{\M}(\interval[0,1];\real)$\@.  Consider the trivial
$\integerp$-net $\map{\phi_{[f]}}{\integerp}{\hat{\M}(\interval[0,1];\real)}$
defined by $\phi_{[f]}(j)=[f]$\@.  Since $\mathcal{F}_\phi=\mathcal{F}_{[f]}$
and since $(\mathcal{F}_\phi,[f])\in\mathscr{L}_\lambda$\@, the
condition~\eqref{pl:ls1} for a limit structure is satisfied.

Let $(\mathcal{F},[f])\in\mathscr{L}_\lambda$ and suppose that
$\mathcal{F}\subset\mathcal{G}$\@.  Then $\mathcal{F}\in\mathscr{F}_{[f]}$
and so $\mathcal{F}\supset\mathcal{F}_\phi$ for some $\integerp$-net $\phi$
that converges pointwise almost everywhere to $[f]$\@.  Therefore, we
immediately have $\mathcal{F}_\phi\subset\mathcal{G}$ and so
$(\mathcal{G},[f])\in\mathscr{L}_\lambda$\@.  This verifies
condition~\eqref{pl:ls2} in the definition of a limit structure.

Finally, let $(\mathcal{F},[f]),(\mathcal{G},[f])\in\mathscr{L}_\lambda$ and
let $\phi$ and $\psi$ be $\integerp$-nets that converge pointwise almost
everywhere to $[f]$ and satisfy $\mathcal{F}_\phi\subset\mathcal{F}$ and
$\mathcal{F}_\psi\subset\mathcal{G}$\@.  Define a $\integerp$-net
$\phi\wedge\psi$ by
\begin{equation*}
\phi\wedge\psi(j)=\begin{cases}\phi(\frac{1}{2}(j+1)),&j\ \textrm{odd},\\
\psi(\frac{1}{2}j),&j\ \textrm{even}.\end{cases}
\end{equation*}
We first claim that $\phi\wedge\psi$ converges pointwise almost everywhere to
$[f]$\@.  Let
\begin{equation*}
A=\asetdef{x\in\interval[0,1]}{\lim_{j\to\infty}\phi(j)(x)\not=f(x)},\quad
B=\asetdef{x\in\interval[0,1]}{\lim_{j\to\infty}\psi(j)(x)\not=f(x)}.
\end{equation*}
If $x\in\interval[0,1]\setminus(A\cup B)$ then
\begin{equation*}
\lim_{j\to\infty}\phi(j)(x)=\lim_{j\to\infty}\psi(j)(x)=f(x).
\end{equation*}
Thus, for $x\in\interval[0,1]\setminus(A\cup B)$ and $\epsilon\in\realp$
there exists $N\in\integerp$ such that
\begin{equation*}
\snorm{f(x)-\phi(j)(x)},\snorm{f(x)-\psi(j)(x)}<\epsilon,\qquad j\ge N.
\end{equation*}
Therefore, for $j\ge2N$ and for $x\in\interval[0,1]\setminus(A\cup B)$ we
have $\snorm{f(x)-\phi\wedge\psi(j)(x)}<\epsilon$ and so
\begin{equation*}
\lim_{j\to\infty}\phi\wedge\psi(j)(x)=f(x),\qquad
x\in\interval[0,1]\setminus(A\cup B).
\end{equation*}
Since $\lambda(A\cup B)=0$ it indeed follows that $\phi\wedge\psi$ converges
pointwise almost everywhere to $[f]$\@.

We next claim that
$\mathcal{F}_{\phi\wedge\psi}\subset\mathcal{F}\cap\mathcal{G}$\@. Indeed,
let $S\in\mathcal{F}_{\phi\wedge\psi}$\@. Then there exists $N\in\integerp$
such that $T_{\phi\wedge\psi}(N)\subset S$\@. Therefore, there exists
$N_\phi,N_\psi\in\integerp$ such that $T_\phi(N_\phi)\subset S$ and
$T_\psi(N_\psi)\subset S$\@. That is,
$S\in\mathcal{F}_\phi\cap\mathcal{F}_\psi\subset\mathcal{F}\cap\mathcal{G}$\@.
This shows that $(\mathcal{F}\cap\mathcal{G},[f])\in\mathscr{L}_\lambda$ and
so shows that condition~\eqref{pl:ls3} in the definition of a limit structure
holds.

Thus we have shown that $\mathscr{L}_\lambda$ is a limit structure.  Let us
show that it is a linear limit structure.  Let
$(\mathcal{F}_1,[f_1]),(\mathcal{F}_2,v_2)\in\mathscr{L}_\lambda$\@.  Thus
there exists $\integer$-nets $\phi_1$ and $\phi_2$ in
$\hat{\M}(\interval[0,1];\real)$ converging pointwise almost everywhere to
$[f_1]$ and $[f_2]$\@, respectively, and such that
$\mathcal{F}_{\phi_1}\subset\mathcal{F}_1$ and
$\mathcal{F}_{\phi_2}\subset\mathcal{F}_2$\@.  Let us denote by
$\ifam{f_{1,j}}_{j\in\integerp}$ and $\ifam{f_{2,j}}_{j\in\integerp}$
sequences in $\M(\interval[0,1];\real)$ such that $[f_{1,j}]=\phi_1(j)$ and
$[f_{2,j}]=\phi_2(j)$ for $j\in\integerp$\@.  Then, as in the proof of
Lemma~\ref{lem:pwaec}\@, there exists a subset $A\subset\interval[0,1]$ of
zero measure such that
\begin{equation*}
\lim_{j\to\infty}f_{j,1}(x)=f_1(x),\quad\lim_{j\to\infty}f_{2,j}(x)=f_2(x),
\qquad x\in\interval[0,1]\setminus A.
\end{equation*}
Thus, for $x\in\interval[0,1]\setminus A$\@,
\begin{equation*}
\lim_{j\to\infty}(f_{1,j}+f_{2,j})(x)=(f_1+f_2)(x).
\end{equation*}
This shows that the $\integerp$-net $\phi_1+\phi_2$ converges pointwise
almost everywhere to $[f_1+f_2]$\@. Since
$\mathcal{F}_{\phi_1+\phi_2}\subset\mathcal{F}_1+\mathcal{F}_2$\@, it follows
that $(\mathcal{F}_1+\mathcal{F}_2,[f_1+f_2])\in\mathscr{L}_\lambda$\@. An
entirely similarly styled argument gives
$(a\mathcal{F},av)\in\mathscr{L}_\lambda$ for
$(\mathcal{F},v)\in\mathscr{L}_\lambda$\@.

We now need to show that $\mathscr{S}(\mathscr{L}_\lambda)$ consists exactly
of the almost everywhere pointwise convergent sequences.  The very definition
of $\mathscr{L}_\lambda$ ensures that if a $\integerp$-net $\phi$ is almost
everywhere pointwise convergent then
$\phi\in\mathscr{S}(\mathscr{L}_\lambda)$\@.  We prove the converse, and so
let $\phi$ be $\mathscr{L}_\lambda$-convergent to $[f]$\@.  Therefore, by
definition of $\mathscr{L}_\lambda$\@, there exists a $\integerp$-net $\psi$
converging pointwise almost everywhere to $[f]$ such that
$\mathcal{F}_\psi\subset\mathcal{F}_\phi$\@.
\begin{prooflemma}
There exists of a subsequence $\psi'$ of $\psi$ such that
$\mathcal{F}_{\psi'}=\mathcal{F}_\phi$\@.
\begin{subproof}
Let $n\in\integerp$ and note that
$T_\psi(n)\in\mathcal{F}_\psi\subset\mathcal{F}_\phi$\@.  Thus there exists
$k\in\integerp$ such that $T_\phi(k)\subset T_\psi(n)$\@.  Then define
\begin{equation*}
k_n=\min\setdef{k\in\integerp}{T_\phi(k)\subset T_\psi(n)},
\end{equation*}
the minimum being well-defined since
\begin{equation*}
k>k'\quad\implies\quad T_\phi(k)\subset T_\phi(k').
\end{equation*}
This uniquely defines, therefore, a sequence $\ifam{k_n}_{n\in\integerp}$\@.
Moreover, if $n_1>n_2$ then $T_\psi(n_2)\subset T_\psi(n_1)$ which implies
that $T_\phi(k_{n_2})\subset T_\psi(n_1)$\@.  Therefore, $k_{n_2}\ge
k_{n_1}$\@, showing that the sequence $\ifam{k_n}_{n\in\integerp}$ is
nondecreasing.

Now define $\map{\theta}{\integerp}{\integerp}$ as follows.  If $j<k_n$ for
every $n\in\integerp$ then define $\theta(j)$ in an arbitrary manner.  If
$j\ge k_1$ then note that $\phi(j)\in T_\phi(k_1)\subset T_\psi(1)$\@.  Thus
there exists (possibly many) $m\in\integerp$ such that $\phi(j)=\psi(m)$\@.
If $j\ge k_n$ for $n\in\integerp$ then there exists (possibly many) $m\ge n$
such that $\phi(j)=\psi(m)$\@.  Thus for any $j\in\integerp$ we can define
$\theta(j)\in\integerp$ such that $\phi(j)=\psi(\theta(j))$ if $j\ge k_1$ and
such that $\theta(j)\ge n$ if $j\ge k_n$\@.

Note that any function $\map{\theta}{\integerp}{\integerp}$ as constructed
above is unbounded.  Therefore, there exists a strictly increasing function
$\map{\rho}{\integerp}{\integerp}$ such that $\image(\rho)=\image(\theta)$\@.
We claim that $\mathcal{F}_\rho=\mathcal{F}_\theta$\@.  First let
$n\in\integerp$ and let $j\ge k_{\rho(n)}$\@.  Then $\theta(j)\ge\rho(n)$\@.
Since $\image(\rho)=\image(\theta)$ there exists $m\in\integerp$ such that
$\rho(m)=\theta(j)\ge\rho(n)$\@.  Since $\rho$ is strictly increasing, $m\ge
n$\@.  Thus $\theta(j)\in T_\rho(n)$ and so $T_\theta(k_{\rho(n)})\subset
T_\rho(n)$\@.  This implies that
$\mathcal{F}_\rho\subset\mathcal{F}_\theta$\@.

Conversely, let $n\in\integerp$ and let $r_n\in\integerp$ be such that
\begin{equation*}
\rho(r_n)>\max\{\theta(1),\dots,\theta(n)\};
\end{equation*}
this is possible since $\rho$ is unbounded.  If $j\ge r_n$ then
\begin{equation*}
\rho(j)\ge\rho(r_n)>\max\{\theta(1),\dots,\theta(n)\}.
\end{equation*}
Since $\image(\rho)=\image(\theta)$ we have $\rho(j)=\theta(m)$ for some
$m\in\integerp$\@.  We must have $m>n$ and so $\rho(j)\in T_\theta(n)$\@.
Thus $T_\rho(r_n)\subset T_\theta(n)$ and so
$\mathcal{F}_\theta\subset\mathcal{F}_\rho$\@.

To arrive at the conclusions of the lemma we first note that, by definition
of $\theta$\@, $\mathcal{F}_\phi=\mathcal{F}_{\psi\scirc\theta}$\@.  We now
define $\psi'=\psi\scirc\rho$ and note that
\begin{equation*}
\mathcal{F}_\phi=\mathcal{F}_{\psi\scirc\theta}=
\psi(\mathcal{F}_\theta)=\psi(\mathcal{F}_\rho)=\mathcal{F}_{\psi\scirc\rho},
\end{equation*}
as desired.
\end{subproof}
\end{prooflemma}

Since a subsequence of an almost everywhere pointwise convergent sequence is
almost everywhere pointwise convergent to the same limit, it follows that
$\psi'$\@, and so $\phi$\@, converges almost everywhere pointwise to $[f]$\@.
\end{proof}
\end{theorem}

The preceding theorem seems to be well-known; see~\cite{RB/HPB:02} where, in
particular, the essential lemma in the proof is given.  Nonetheless, we have
never seen the ingredients of the proof laid out clearly in one place, so the
result is worth recording.

Let us record a characterisation of $\mathscr{L}_\lambda$-bounded subsets of
$\hat{\M}(\interval[0,1];\real)$\@.
\begin{proposition}\label{prop:hatM-bounded}
A subset\/ $B\subset\hat{\M}(\interval[0,1];\real)$ is\/
$\mathscr{L}_\lambda$-bounded if and only if there exists a
nonnegative-valued\/ $g\in\M(\interval[0,1];\real)$ such that
\begin{equation*}
B\subset\setdef{[f]\in\hat{\M}(\interval[0,1];\real)}
{\snorm{f(x)}\le g(x)\ \textrm{for almost every}\ x\in\interval[0,1]}.
\end{equation*}
\begin{proof}
We first observe that the condition that $\snorm{f(x)}\le g(x)$ for almost
every $x\in\interval[0,1]$ is independent of the choice of representative $f$
from the equivalence class $[f]$\@.

Suppose that there exists a nonnegative-valued\/
$g\in\M(\interval[0,1];\real)$ such that, if\/ $[f]\in B$\@, then\/
$\snorm{f(x)}\le g(x)$ for almost every\/ $x\in\interval[0,1]$\@.  Let
$\ifam{[f_j]}_{j\in\integerp}$ be a sequence in $B$ and let
$\ifam{a_j}_{j\in\integerp}$ be a sequence in $\real$ converging to zero.
For $j\in\integerp$ define
\begin{equation*}
A_j=\setdef{x\in\interval[0,1]}{\snorm{f_j(x)}\le g(x)}.
\end{equation*}
Note that if $x\in\interval[0,1]\setminus(\cup_{j\in\integerp}A_j)$ then
\begin{equation*}
\lim_{j\to\infty}\snorm{a_jf_j(x)}\le\lim_{j\to\infty}\snorm{a_j}g(x)=0.
\end{equation*}
Since $\lambda(\cup_{j\in\integerp}A_j)=0$ this implies that the sequence
$\ifam{a_j[f_j]}_{j\in\integerp}$ is $\mathscr{L}_\lambda$-convergent to
zero.  One may show that this argument is independent of the choice of
representatives $f_j$ from the equivalence classes $[f_j]$\@,
$j\in\integerp$\@.

Conversely, suppose that there exists no nonnegative-valued function
$g\in\M(\interval[0,1];\real)$ such that, for every $[f]\in B$\@,
$\snorm{f(x)}\le g(x)$ for almost every $x\in\interval[0,1]$\@.  This means
that there exists a set $E\subset\interval[0,1]$ of positive measure such
that, for any $M\in\realp$\@, there exists $[f]\in B$ such that
$\snorm{f(x)}>M$ for almost every $x\in E$\@.  Let
$\ifam{a_j}_{j\in\integerp}$ be a sequence in $\real$ converging to $0$ and
such that $a_j\not=0$ for every $j\in\integerp$\@.  Then let
$\ifam{[f_j]}_{j\in\integerp}$ be a sequence in $B$ such that
$\snorm{f_j(x)}>\asnorm{a_j^{-1}}$ for almost every $x\in E$ and for
every $j\in\integerp$\@.  Define
\begin{equation*}
A_j=\setdef{x\in E}{\snorm{f_j(x)}>\asnorm{a_j^{-1}}}.
\end{equation*}
If $x\in E\setminus(\cup_{j\in\integerp}A_j)$ then $\snorm{a_jf_j(x)}>1$ for
every $j\in\integerp$\@.  Since
$\lambda(E\setminus(\cup_{j\in\integerp}A_j))>0$ it follows that
$\ifam{a_j[f_j]}_{j\in\integerp}$ cannot $\mathscr{L}_\lambda$-converge to
zero, and so $B$ is not $\mathscr{L}_\lambda$-bounded.
\end{proof}
\end{proposition}

\section{Two topological distinctions between the Riemann and Lebesgue
theories of integration}

In this section we give topological characterisations of the differences
between the Riemann and Lebesgue theories. In
Section~\ref{subsec:riemann-ccft} we also explicitly see how these
distinctions lead to a deficiency in the Fourier transform theory using the
Riemann integral.

\subsection{The Dominated Convergence Theorems}

Both the Lebesgue and Riemann theories of integration possess a Dominated
Convergence Theorem.  This gives us two versions of the Dominated Convergence
Theorem that we can compare and contrast.  Moreover, there are also
``pointwise convergent'' and ``almost everywhere pointwise convergent''
versions of both theorems.  Typically, the ``pointwise convergent'' version
is stated for the Riemann integral\footnote{Since the ``almost everywhere
pointwise convergent'' version actually requires the Lebesgue theory of
integration.}  and the ``almost everywhere pointwise convergent'' version is
stated for the Lebesgue integral.  However, both versions are valid for both
integrals, so this gives, in actuality, four theorems to compare and
contrast.  What we do here is state both versions of the Dominated Convergence
Theorem for the Lebesgue integral using topological and limit structures, and
we give counterexamples illustrating why these statements are not valid for
the Riemann integral.

Let us first state the various Dominated Convergence Theorems in their usual
form.  The Dominated Convergence Theorem\textemdash{}including the pointwise
convergent and almost everywhere pointwise convergent
statements\textemdash{}for the Riemann integral is the following.
\begin{theorem}
Let\/ $\ifam{f_j}_{j\in\integerp}$ be a sequence of\/ $\real$-valued
functions on\/ $\interval[0,1]$ satisfying the following conditions:
\begin{compactenum}[(i)]
\item $f_j\in\R^1(\interval[0,1];\real)$ for each\/ $j\in\integerp$\@;
\item there exists a nonnegative-valued\/ $g\in\R^1(\interval[0,1];\real)$
such that\/ $\snorm{f_j(x)}\le g(x)$ for every (\resp~almost every)\/
$x\in\interval[0,1]$ and for every\/ $j\in\integerp$\@;
\item the limit\/ $\lim_{j\to\infty}f_j(x)$ exists for every (\resp~almost
every)\/ $x\in\interval[0,1]$\@;
\item the function\/ $\map{f}{\interval[0,1]}{\real}$ defined by
$f(x)=\lim_{j\to\infty}f_j(x)$ is in\/ $\R^1(\interval[0,1];\real)$
(\resp~there exists\/ $f\in\R^1(\interval[0,1];\real)$ such that\/
$\lim_{j\to\infty}f_j(x)=f(x)$ for almost every\/ $x\in\interval[0,1]$).
\end{compactenum}
Then
\begin{equation*}
\lim_{j\to\infty}\int_0^1f_j(x)\,\d{x}=\int_0^1f(x)\,\d{x}.
\end{equation*}
\end{theorem}

For the Lebesgue integral we have the following Dominated Convergence
Theorem(s).
\begin{theorem}
Let\/ $\ifam{f_j}_{j\in\integerp}$ be a sequence of\/ $\real$-valued
functions on\/ $\interval[0,1]$ satisfying the following conditions:
\begin{compactenum}[(i)]
\item $f_j$ is measurable for each\/ $j\in\integerp$\@;
\item there exists a nonnegative-valued\/ $g\in\L^1(\interval[0,1];\real)$
such that\/ $\snorm{f_j(x)}\le g(x)$ for every (\resp~almost every)\/
$x\in\interval[0,1]$ and for every\/ $j\in\integerp$\@;
\item the limit\/ $\lim_{j\to\infty}f_j(x)$ exists for every (\resp~almost
every)\/ $x\in\interval[0,1]$\@.
\end{compactenum}
Then the function\/ $\map{f}{\interval[0,1]}{\real}$ defined by
\begin{equation*}
f(x)=\begin{cases}\lim_{j\to\infty}f_j(x),&\textrm{the limit exists},\\
0,&\textrm{otherwise}\end{cases}
\end{equation*}
and each of the functions\/ $f_j$\@, $j\in\integerp$\@, are in\/
$\L^1(\interval[0,1];\real)$ and
\begin{equation*}
\lim_{j\to\infty}\int_If_j\,\d{\lebmes}=\int_If\,\d{\lebmes}.
\end{equation*}
\end{theorem}

Our statements make it clear that there is one real difference between the
Riemann and Lebesgue theories:~the condition of integrability of the limit
function $f$ is an \emph{hypothesis} in the Riemann theory but a
\emph{conclusion} in the Lebesgue theory. This distinction is crucial and
explains why the Lebesgue theory is more powerful than the Riemann theory.
Moreover, the structure we introduced in Section~\ref{sec:LR-topologies}
allows for an elegant expression of this distinction. The result which
follows is simply a rephrasing of the Dominated Convergence Theorem(s) for
the Lebesgue integral, and follows from that theorem, along with
Theorem~\ref{the:hatM-convergence} and Proposition~\ref{prop:hatM-bounded}\@.
\begin{theorem}\label{the:DCT} The following statements hold:
\begin{compactenum}[(i)]
\item \label{pl:DCT1} $\C_p$-bounded subsets of\/
$\L^1(\interval[0,1];\real)$ are\/ $\C_p$-sequentially closed;
\item \label{pl:DCT2} $\mathscr{L}_\lambda$-bounded subsets of\/
$\hat{\L}^1(\interval[0,1];\real)$ are\/ $\mathscr{L}_\lambda$-sequentially
closed.
\end{compactenum}
\end{theorem}

The necessity of the weaker conclusions for the Dominated Convergence
Theorem(s) for the Riemann integral is illustrated by the following examples.
First we show why part~\eqref{pl:DCT1} of Theorem~\ref{the:DCT} does not hold
for the Riemann integral.
\begin{example}\label{eg:R1!Cp-closed}
By means of an example, we show that there are $\C_p$-bounded subsets of the
seminormed vector space $\R^1(\interval[0,1];\real)$ that are not
sequentially closed in the topology of $\C_p(\interval[0,1];\real)$\@.  Let
us denote
\begin{equation*}
B=\asetdef{f\in\R^1(\interval[0,1];\real)}{\snorm{f(x)}\le1},
\end{equation*}
noting by Proposition~\ref{prop:Cp-bounded} that $B$ is $\C_p$-bounded.  Let
$\ifam{q_j}_{j\in\integerp}$ be an enumeration of the rational numbers in
$\interval[0,1]$ and define a sequence $\ifam{F_k}_{k\in\integerp}$ in
$\R^1(\interval[0,1];\real)$ by
\begin{equation*}
F_k(x)=\begin{cases}1,&x\in\{q_1,\dots,q_k\},\\
0,&\textrm{otherwise}.\end{cases}
\end{equation*}
The sequence converges in $\C_p(\interval[0,1];\real)$ to the characteristic
function of $\rational\cap\interval[0,1]$\@; let us denote this function by
$F$\@.  This limit function is not Riemann integrable and so not in
$\R^1(\interval[0,1];\real)$\@.  Thus $B$ is not $\C_p$-sequentially
closed.\oprocend
\end{example}

Next we show why part~\eqref{pl:DCT2} of Theorem~\ref{the:DCT} does not hold
for the Riemann integral.
\begin{example}\label{eg:R1!hatMp-closed}
We give an example that shows that $\mathscr{L}_\lambda$-bounded subsets of
the normed vector space $\hat{\R}^1(\interval[0,1];\real)$ are not
$\mathscr{L}_\lambda$-sequentially closed.  We first remark that the
construction of Example~\ref{eg:R1!Cp-closed}\@, projected to
$\hat{\R}^1(\interval[0,1];\real)$\@, does not suffice because $[F]$ is equal
to the equivalence class of the zero function which \emph{is} Riemann
integrable, even though $F$ is not,~\cf~the statements following
Example~3.117 in~\cite{DSK:04}\@.  The fact that $[F]$ contains functions
that are Riemann integrable and functions that are not Riemann integrable is
a reflection of the fact that $\R_0(\interval[0,1];\real)$ is not
sequentially closed.  This is a phenomenon of interest, but it is not what is
of interest here.

The construction we use is the following.  Let $\ifam{q_j}_{j\in\integerp}$
be an enumeration of the rational numbers in $\interval[0,1]$\@.  Let
$\ell\in\interval(0,1)$ and for $j\in\integerp$ define
\begin{equation*}
I_j=\interval[0,1]\cap
\interval({q_j-\tfrac{\ell}{2^{j+1}}},{q_j+\tfrac{\ell}{2^{j+1}}})
\end{equation*}
to be the interval of length $\frac{\ell}{2^j}$ centred at $q_j$\@.  Then
define $A_k=\cup_{j=1}^kI_j$\@, $k\in\integerp$\@, and
$A=\cup_{j\in\integerp}A_j$\@.  Also define $G_k=\chi_{A_k}$\@,
$k\in\integerp$\@, and $G=\chi_A$ be the characteristic functions of $A_k$
and $A$\@, respectively.  Note that $A_k$ is a union of a finite number of
intervals and so $G_k$ is Riemann integrable for each $k\in\integerp$\@.
However, we claim that $G$ is not Riemann integrable.  Indeed, the
characteristic function of a set is Riemann integrable if and only the
boundary of the set has measure zero; this is a direct consequence of
Lebesgue's theorem stating that a function is Riemann integrable if and only
if its set of discontinuities has measure
zero~\cite[Theorem~2.5.4]{DLC:13}\@.  Note that since
$\closure(\rational\cap\interval[0,1])=\interval[0,1]$ we have
\begin{equation*}
\interval[0,1]=\closure(A)=A\cup\boundary(A).
\end{equation*}
Thus
\begin{equation*}
\lebmes(\interval[0,1])\le\lebmes(A)+\lebmes(\boundary(A)).
\end{equation*}
Since
\begin{equation*}
\lebmes(A)\le\sum_{j=1}^\infty\lebmes(I_j)\le\ell,
\end{equation*}
it follows that $\lebmes(\boundary(A))\ge1-\ell\in\realp$\@.  Thus $G$ is
not Riemann integrable, as claimed.

It is clear that $\ifam{G_k}_{k\in\integerp}$ is $\C_p$-convergent to $G$\@.
Therefore, by Theorem~\ref{the:hatM-convergence} it follows that
$\ifam{[G_k]}_{k\in\integerp}$ is $\mathscr{L}_\lambda$-convergent to
$[G]$\@.  We claim that $[G]\not\in\hat{\R}^1(\interval[0,1];\real)$\@.  This
requires that we show that if $g\in\C_p(\interval[0,1];\real)$ satisfies
$[g]=[G]$\@, then $g$ is not Riemann integrable.  To show this, it suffices
to show that $g$ is discontinuous on a set of positive measure.  We shall
show that $g$ is discontinuous on the set $g^{-1}(0)\cap\boundary(A)$\@.
Indeed, let $x\in g^{-1}(0)\cap\boundary(A)$\@.  Then, for any
$\epsilon\in\realp$ we have $\interval({x-\epsilon},{x+\epsilon})\cap
A\not=\emptyset$ since $x\in\boundary(A)$\@.  Since
$\interval({x-\epsilon},{x+\epsilon})\cap A$ is a nonempty open set, it has
positive measure.  Therefore, since $G$ and $g$ agree almost everywhere,
there exists $y\in\interval({x-\epsilon},{x+\epsilon})\cap A$ such that
$g(y)=1$\@.  Since this holds for every $\epsilon\in\realp$ and since
$g(x)=0$\@, it follows that $g$ is discontinuous at $x$\@.  Finally, it
suffices to show that $g^{-1}(0)\cap\boundary(A)$ has positive measure.  But
this follows since $\boundary(A)=G^{-1}(0)$ has positive measure and since
$G$ and $g$ agree almost everywhere.

To complete the example, we note that the sequence
$\ifam{[G_k]}_{j\in\integerp}$ is in the set
\begin{equation*}
B=\asetdef{[f]\in\hat{\R}^1(\interval[0,1];\real)}
{\snorm{f(x)}\le 1\ \textrm{for almost every}\ x\in\interval[0,1]},
\end{equation*}
which is $\mathscr{L}_\lambda$-bounded by
Proposition~\ref{prop:hatM-bounded}\@.  The example shows that this
$\mathscr{L}_\lambda$-bounded subset of $\hat{\R}^1(\interval[0,1];\real)$ is
not $\mathscr{L}_\lambda$-sequentially closed.\oprocend
\end{example}

\subsection{Completeness of spaces of integrable functions}

In Section~\ref{subsec:nvs-integrable} we constructed the two normed vector
spaces $\hat{\R}^1(\interval[0,1];\real)$ and
$\hat{\L}^1(\interval[0,1];\real)$\@.  Using the Dominated Convergence
Theorem, one proves the following well-known and important
result~\cite[Theorem~3.4.1]{DLC:13}\@.
\begin{theorem}
$(\hat{\L}^1(\interval[0,1];\real),\dnorm{\cdot}_1)$ is a Banach space.
\end{theorem}

It is generally understood that
$(\hat{\R}^1(\interval[0,1];\real),\dnorm{\cdot}_1)$ is \emph{not} a Banach
space.  However, we have not seen this demonstrated in a sufficiently
compelling manner, so the following example will hopefully be interesting in
this respect.
\begin{example}\label{eg:Riemann!Cauchy}
Let us consider the sequence of functions $\ifam{G_k}_{k\in\integerp}$ in
$\R^1(\interval[0,1];\real)$ constructed in
Example~\ref{eg:R1!hatMp-closed}\@.  We also use the pointwise limit function
$G$ defined in that same example.  We shall freely borrow the notation
introduced in this example.

We claim that the sequence $\ifam{[G_k]}_{k\in\integerp}$ is Cauchy in
$\hat{\R}^1(\interval[0,1];\real)$\@.  Let $\epsilon\in\realp$\@.  Note that
$\sum_{j=1}^\infty\lebmes(I_j)\le\ell$\@.  This implies that there exists
$N\in\integerp$ such that $\sum_{j=k+1}^m\lebmes(I_j)<\epsilon$ for all
$k,m\ge N$\@.  Now note that for $k,m\in\integerp$ with $m>k$\@, the
functions $G_k$ and $G_m$ agree except on a subset of $I_{k+1}\cup\dots\cup
I_m$\@.  On this subset, $G_m$ has value $1$ and $G_k$ has value $0$\@.  Thus
\begin{equation*}
\int_0^1\snorm{G_m(x)-G_k(x)}\,\d{x}\le\lebmes(I_{k+1}\cup\dots\cup I_m)\le
\sum_{j=k+1}^m\lebmes(I_j).
\end{equation*}
Thus we can choose $N\in\integerp$ sufficiently large that
$\dnorm{G_m-G_k}_1<\epsilon$ for $k,m\ge N$\@.  Thus the sequence
$\ifam{[G_k]}_{k\in\integerp}$ is Cauchy, as claimed.

We next show that the sequence $\ifam{[G_k]}_{k\in\integerp}$ converges to
$[G]$ in $\hat{\L}^1(\interval[0,1];\real)$\@.  In
Example~\ref{eg:R1!hatMp-closed} we showed that
$\ifam{[G_k]}_{k\in\integerp}$ is $\mathscr{L}_\lambda$-convergent to $[G]$).
Since the sequence $\ifam{[G-G_k]}_{k\in\integerp}$ is in the subset
\begin{equation*}
\setdef{[f]\in\hat{\L}^1(\interval[0,1];\real)}
{\snorm{f(x)}\le1\ \textrm{for almost every}\ x\in\interval[0,1]},
\end{equation*}
and since this subset is $\mathscr{L}_\lambda$-bounded by
Proposition~\ref{prop:hatM-bounded}\@, it follows from
Theorem~\pldblref{the:DCT}{pl:DCT2} that
\begin{equation*}
\lim_{k\to\infty}\dnorm{G-G_k}_1=
\int_I\lim_{k\to\infty}\snorm{G-G_k}\,\d{\lambda}=0.
\end{equation*}
This gives us the desired convergence of $\ifam{[G_k]}_{k\in\integerp}$ to
$[G]$ in $\hat{\L}^1(\interval[0,1];\real)$\@.  However, in
Example~\ref{eg:R1!hatMp-closed} we showed that
$[G]\not\in\hat{\R}^1(\interval[0,1];\real)$\@.  Thus the Cauchy sequence
$\ifam{[G_k]}_{k\in\integerp}$ in $\hat{\R}^1(\interval[0,1];\real)$ is not
convergent in $\hat{\R}(\interval[0,1];\real)$\@, giving the desired
incompleteness of
$(\hat{\R}(\interval[0,1];\real),\dnorm{\cdot}_1)$\@.\oprocend
\end{example}

In \cite[Example~3.117]{DSK:04} the sequence $\ifam{f_j}_{j\in\integerp}$
defined by
\begin{equation*}
f_j(x)=\begin{cases}0,&x\in\interval[0,{\frac{1}{j}}],\\
x^{-1/2},&x\in\interval({\frac{1}{j}},1]\end{cases}
\end{equation*}
is shown to be Cauchy in $\hat{\R}^1(\interval[0,1];\real)$\@, but not
convergent in $\hat{\R}^1(\interval[0,1];\real)$\@.  This sequence, however,
is not as interesting as that in our preceding example since the limit
function $f\in\hat{\L}^1(\interval[0,1];\real)$ in \citeauthor{DSK:04}'s
example \emph{is} Riemann integrable using the usual rule for defining the
improper Riemann integral for unbounded functions.  In the construction used
in Example~\ref{eg:Riemann!Cauchy}\@, the limit function in
$\L^1(\interval[0,1];\real)$ is not Riemann integrable in the sense of
bounded functions defined on compact intervals,~\ie~in the sense of the usual
construction involving approximation above and below by step functions.

In \cite{WD/HM/RJN/EvW:88} a convergence structure is introduced on the set
of Riemann integrable functions that is sequentially complete.  The idea in
this work is to additionally constrain convergence in
$\hat{\L}^1(\interval[0,1];\real)$ in such a way that Riemann integrability
is preserved by limits.

\subsection{The $\L^2$-Fourier transform for the Riemann
integral}\label{subsec:riemann-ccft}

In order to illustrate why it is important that spaces of integrable
functions be complete, we consider the theory of the $\L^2$-Fourier transform
restricted to square Riemann integrable functions.  Let us first recall the
essential elements of the theory.

For $p\in\interval[1,\infty)$ we denote by $\L^p(\real;\complex)$ the set of
$\complex$-valued functions $f$ defined on $\real$ which satisfy
\begin{equation*}
\int_{\real}\snorm{f}^p\,\d{\lambda}<\infty,
\end{equation*}
where we now denote by $\lambda$ the Lebesgue measure on $\real$\@.  We let
\begin{equation*}
\L_0(\real;\complex)=\asetdef{f\in\L^p(\real;\complex)}
{\int_{\real}\snorm{f}^p\,\d{\lambda}=0}
\end{equation*}
and denote
\begin{equation*}
\hat{\L}^p(\real;\complex)=\L^p(\real;\complex)/\L_0(\real;\complex).
\end{equation*}
As we have done previously, we denote $[f]=f+\L_0(\real;\complex)$\@.  If we
define
\begin{equation*}
\dnorm{[f]}_p=\left(\int_{\real}\snorm{f}^p\,\d{\lambda}\right)^{1/p}
\end{equation*}
then one shows that $(\hat{\L}^p(\real;\complex),\dnorm{\cdot}_p)$ is a
Banach space~\cite[Theorem~3.4.1]{DLC:13}\@.

For $a\in\complex$ let us denote $\map{\exp_a}{\real}{\complex}$ by
$\exp_a(x)=\eul^{ax}$\@.  For $[f]\in\hat{\L}^1(\real;\complex)$ one defines
$\map{\mathscr{F}_1([f])}{\real}{\complex}$ by
\begin{equation*}
\mathscr{F}_1([f])(\omega)=\int_{\real}f\exp_{-2\pi\imag\omega}\,\d{\lambda}
\end{equation*}
the \defn{$\L^1$-Fourier transform} of $[f]\in\hat{\L}^1(\real;\complex)$\@.
If we define $\C^0_{0,\textup{uc}}(\real;\complex)$ to be the set of
uniformly continuous $\complex$-valued functions $f$ on $\real$ that satisfy
$\lim_{\snorm{x}\to\infty}\snorm{f(x)}=0$\@, then
$(\C^0_{0,\textup{uc}}(\real;\complex),\dnorm{\cdot}_\infty)$ is a Banach
space with
\begin{equation*}
\dnorm{f}_\infty=\sup\setdef{\snorm{f(x)}}{x\in\real}.
\end{equation*}
Moreover, $\mathscr{F}_1([f])\in\C^0_{0,\textup{uc}}(\real;\complex)$ and the
linear map $\map{\mathscr{F}_1}{\hat{\L}^1(\real;\complex)}
{\C^0_{0,\textup{uc}}(\real;\complex)}$ is
continuous~\cite[Theorem~17.1.3]{CG/PW:99}\@.

For $[f]\in\hat{\L}^2(\real;\complex)$ the preceding construction cannot be
applied verbatim since $\hat{\L}^2(\real;\complex)$ is not a subspace of
$\hat{\L}^1(\real;\complex)$\@.  However, one can make an adaptation as
follows~\cite[Lesson~22]{CG/PW:99}\@.  One shows that
$\hat{\L}^1(\real;\complex)\cap\hat{\L}^2(\real;\complex)$ is dense in
$\hat{\L}^2(\real;\complex)$\@.  One can do this explicitly by defining, for
$[f]\in\hat{\L}^2(\real;\complex)$\@, a sequence
$\ifam{[f_j]}_{j\in\integerp}$ in
$\hat{\L}^1(\real;\complex)\cap\hat{\L}^2(\real;\complex)$ by
\begin{equation*}
f_j(x)=\begin{cases}f(x),&x\in\interval[-j,j],\\
0,&\textrm{otherwise},\end{cases}
\end{equation*}
and showing, using the Cauchy\textendash{}Bunyakovsky\textendash{}Schwarz
inequality, that this sequence converges in $\hat{\L}^2(\real;\complex)$ to
$[f]$\@.  Moreover, one can show that the sequence
$\ifam{\mathscr{F}_1([f_j])}_{j\in\integerp}$ is a Cauchy sequence
$\hat{\L}^2(\real;\complex)$ and so converges to some element of
$\hat{\L}^2(\real;\complex)$ that we denote by $\mathscr{F}_2([f])$\@, the
\defn{$\hat{\L}^2$-Fourier transform} of
$[f]\in\hat{\L}^2(\real;\complex)$\@.  The map
$\map{\mathscr{F}_2}{\hat{\L}^2(\real;\complex)}{\hat{\L}^2(\real;\complex)}$
so defined is, moreover, a Hilbert space isomorphism.  The inverse has the
property that
\begin{equation*}
\mathscr{F}_2^{-1}([f])(x)=\int_{\real}f\exp_{2\pi\imag x}\,\d{\lambda}
\end{equation*}
for almost every $x\in\real$\@, where the same constructions leading to the
definition of $\mathscr{F}_2$ for functions in $\hat{\L}^2(\real;\complex)$
are applied.

Let us see that $\mathscr{F}_2$ restricted to the subspace of square Riemann
integrable functions is problematic.  We denote by $\R^p(\real;\complex)$ the
collection of functions $\map{f}{\real}{\complex}$ which satisfy
\begin{equation*}
\int_{-\infty}^\infty\snorm{f(x)}^p\,\d{x}<\infty,
\end{equation*}
where we use, as above, the Riemann integral for possibly unbounded functions
defined on unbounded domains~\cite[Section~8.5]{JEM/MJH:93}\@.  We also
define
\begin{equation*}
\R_0(\real;\complex)=\asetdef{f\in\R^p(\real;\complex)}
{\int_{-\infty}^\infty\snorm{f(x)}^p\,\d{x}=0}
\end{equation*}
and denote
\begin{equation*}
\hat{\R}^p(\real;\complex)=\R^p(\real;\complex)/\R_0(\real;\complex).
\end{equation*}
As we have done previously, we denote $[f]=f+\R_0(\real;\complex)$\@.  If we
define
\begin{equation*}
\dnorm{[f]}_p=\left(\int_{-\infty}^\infty\snorm{f(x)}^p\,\d{x}\right)^{1/p}
\end{equation*}
then $(\hat{\R}^p(\real;\complex),\dnorm{\cdot}_p)$ is a normed vector space.
It is not a Banach space since the example of Example~\ref{eg:Riemann!Cauchy}
can be extended to $\hat{\R}^p(\real;\complex)$ by taking all functions to be
zero outside the interval $\interval[0,1]$\@.

Let us show that $\mathscr{F}_2|\hat{\R}^2(\real;\complex)$ does not take
values in $\hat{\R}^2(\real;\complex)$\@, and thus show that the
``$\hat{\R}^2$-Fourier transform'' is not well-defined.  We denote by $G$ the
function defined in Example~\ref{eg:R1!hatMp-closed}\@, but now extended to
be defined on $\real$ by taking it to be zero off $\interval[0,1]$\@.  We
have $G\in\L^1(\real;\complex)\cap\L^2(\real;\complex)$ since $G$ is bounded
and measurable with compact support.  Now define $\map{F}{\real}{\complex}$
by
\begin{equation*}
F(x)=\int_{\real}G\exp_{2\pi\imag x}\,\d{\lambda};
\end{equation*}
thus $F$ is the inverse Fourier transform of $G$\@.  Since
$G\in\L^1(\real;\complex)$ it follows that
$F\in\C^0_{0,\textup{uc}}(\real;\complex)$\@.  Therefore, $F|\interval[-R,R]$
is continuous and bounded, and hence Riemann integrable for every
$R\in\realp$\@.  Since $G\in\L^2(\real;\complex)$ we have
$F\in\L^2(\real;\complex)$ which implies that
\begin{equation*}
\int_{-R}^R\snorm{F(x)}^2\,\d{x}=
\int_{\interval[-R,R]}\snorm{F}^2\,\d{\lambda}\le
\int_{\real}\snorm{F}^2\,\d{\lambda},\qquad R\in\realp.
\end{equation*}
Thus the limit
\begin{equation*}
\lim_{R\to\infty}\int_{-R}^R\snorm{F(x)}^2\,\d{x}
\end{equation*}
exists.  This is exactly the condition for Riemann integrability of $F$ as a
function on an unbounded domain~\cite[Section~8.5]{JEM/MJH:93}\@.  Now, since
$[F]=\mathscr{F}_2^{-1}([G])$ by definition, we have
$\mathscr{F}_2([F])=[G]$\@.  In Example~\ref{eg:R1!hatMp-closed} we showed
that $[G]|\interval[0,1]\not\in\hat{\R}^1(\interval[0,1];\complex)$\@.  From
this we conclude that $[G]\not\in\hat{\R}^1(\real;\complex)$ and, since
$\snorm{G}^2=G$\@, $[G]\not\in\hat{\R}^2(\real;\complex)$\@.  Thus
$\mathscr{F}_2(\hat{\R}^2(\real;\complex))\not\subset
\hat{\R}^2(\real;\complex)$\@, as it was desired to show.

\printbibliography
\end{document}